\documentclass[11pt]{article}
\usepackage{latexsym, amsfonts,mathrsfs, amsmath, amssymb, amscd, epsfig}

\newtheorem{thm}{Theorem}

\newtheorem{lem}[thm]{Lemma}

\newtheorem{rmk}{Remark}

\newenvironment{pf}{{\noindent \it \bf Proof:}}{{\hfill$\Box$}\\}

\begin{document}

\def\del{\partial}
\def\DOT{\!\cdot\!}
\def\dt{{\Delta t}}
\def\dx{{\Delta x}}
\def\dy{{\Delta y}}
\def\dz{{\Delta z}}
\def\div{\text{div}}
\def\u{{\mbox{\boldmath $u$}}}
\def\g{{\mbox{\boldmath $g$}}}
\def\bu{{\mbox{\boldmath $u$}}}
\def\x{{\mbox{\boldmath $x$}}}
\def\n{{\mbox{\boldmath $n$}}}
\def\bn{{\mbox{\boldmath $n$}}}
\def\veceps{\mbox{\boldmath $\epsilon$}}
\def\vecdel{\mbox{\boldmath $\delta$}}
\def\vecphi{\mbox{\boldmath $\phi$}}
\def\vecpsi{\mbox{\boldmath $\psi$}}
\def\ds{\mbox{d \boldmath $\gamma$}}
\def\vecomi{\mbox{\boldmath $\theta$}}
\def\eps{\varepsilon}
\def\oneb{{1\mbox \tiny b}}
\def\twob{{2\mbox \tiny b}}
\def\disp{\displaystyle}
\def\dspace{\displaystyle \vspace{.15in}}
\def\half{\textstyle \frac12}
\def\Circ{{\footnotesize$\bigcirc$}}
\def\Triangle{{\small$\triangle$}}
\def\Bullet{{\large$\bullet$}}
\def\exac{{\mbox \tiny e}}
\newcommand{\Dtil}{\widetilde{D}}

\def\TIME{\!\times\!}

\title{Global Subsonic and Subsonic-Sonic Flows\\
through Infinitely Long Axially Symmetric Nozzles}

\author {Chunjing \, Xie\footnote{Department of mathematics, University of Michigan,
Ann Arbor, E-mail address: cjxie@umich.edu}\,\, and \,\, Zhouping \,
Xin\footnote{The Institute of Mathematical Sciences and department
of mathematics, The Chinese University of Hong Kong, E-mail address:
zpxin@ims.cuhk.edu.hk}}

\date{}
\maketitle

{\bf Abstract:} In this paper, we establish existence of global
subsonic and subsonic-sonic flows through infinitely long axially
symmetric nozzles by combining variational method, various elliptic
estimates and a compensated compactness method. More precisely, it
is shown that there exist global subsonic flows in nozzles for
incoming mass flux less than a critical value; moreover, uniformly
subsonic flows always approach to uniform flows at far fields when
nozzle boundaries tend to be flat at far fields, and flow angles for
axially symmetric flows are uniformly bounded away from $\pi/2$;
finally, when the incoming mass flux tends to the critical value,
subsonic-sonic flows exist globally in nozzles in the weak sense by
using angle estimate in conjunction with a compensated compactness
framework.


 \noindent

\bigskip

\section{Introduction and Main Results}

This is a continuation to \cite{XX} on the study for subsonic and
subsonic-sonic flows in multidimensional nozzles. Two dimensional
subsonic and subsonic-sonic flows through infinitely long nozzles
were studied in detail in \cite{XX}. In particular, global smooth
subsonic flows and their properties have been established for
incoming mass fluxes less than the critical value, while existence
of subsonic-sonic flows is proved for the critical mass flux.
However, many of arguments are special to 2-dimensional flows, and
it seems difficult to generalize them to the more realistic
3-dimensional flows.

In this paper, we would like to investigate the 3-D flows in
nozzles which are infinitely long and axially symmetric.

As far as axially symmetric flows are concerned, one should note
the significant result due to Gilbarg, \cite{G}, where he showed
that if an axially symmetric subsonic nozzle flow approximates to
uniform flows at far fields, then the flow speed on the boundary
is monotone increasing with respect to the incoming mass flux by a
comparison principle, however, existence of such flows is not
known. For free boundary problems, in \cite{ACF}, Alt, Caffarelli
and Friedman gave a complete study for flows with jet and
cavitation by variational methods.

Let us start with 3-D isentropic compressible Euler equations,
\begin{equation}\label{Euler}
\left\{
\begin{array}{ll}
(\rho u)_x+(\rho v)_y+(\rho w)_z=0,\\
(\rho u^2)_x+(\rho u v)_y+(\rho u w)_z+p(\rho)_x=0,\\
(\rho uv)_x+(\rho  v^2)_y+(\rho v w)_z+p(\rho)_y=0,\\
(\rho uw)_x+(\rho  vw)_y+(\rho  w^2)_z+p(\rho)_z=0,
\end{array}
\right.
\end{equation}
where $\rho$ is the density, $(u,v,w)$ is the velocity, and
$p=p(\rho)$ denotes the pressure. In general, it is assumed that
$p'(\rho)>0$ for $\rho>0$ and $p''(\rho)\geq 0$, where
$c(\rho)=\sqrt{p'(\rho)}$ is called the sound speed. Important
examples include polytropic gases and isothermal gases, for
polytropic gases, $p=A\rho^{\gamma}$ where $A$ is a constant and
$\gamma$ is the adiabatic constant with $\gamma>1$; and for
isothermal gases, $p=c^2\rho$ with constant sound speed $c$.

Suppose that the flow is also irrotational, i.e.\cite{CF},
\begin{equation}\label{3Dirrotational}
u_y-v_x=0,\,\,u_z-w_x,\,\,v_z-w_y=0.
\end{equation}
Then it follows from (\ref{Euler}) and (\ref{3Dirrotational}) that
the flow satisfies Bernoulli's law
\begin{equation}\label{3DBernoulli}
\frac{q^2}{2}+\int^{\rho}\frac{p'(\rho)}{\rho}d\rho=C,
\end{equation}
where $q=\sqrt{u^2+v^2+w^2}$, and $C$ is a constant depending on
the flow. There are some basic facts about irrotational isentropic
steady flows, see \cite{CF}, which are consequences of Bernoulli's
law (\ref{3DBernoulli}). First, $\rho$ is a decreasing function of
$q$, attains its maximum at $q=0$. Second, there is a critical
speed $q_c$ such that $q<c$(subsonic) if and only if $q<q_c$.
Finally, $\rho q$ is a nonnegative function of $q$, for $q \geq
0$, which is increasing for $q\in(0,q_c)$ and decreasing for
$q\geq q_c$, and vanishes at $q=0$. so $\rho q$ attains its
maximum at $q=q_c$, therefore, that the flow is subsonic is
equivalent to $\rho q<\rho_c q_c$ and $\rho>\rho_c$. Therefore, we
can nondimensionalize the flow as in \cite{Bers1,XX}, such that
$q_{cr}=1$, $\rho_{cr}=1$, and then Bernoulli's law
(\ref{3DBernoulli}) reduces to
\begin{equation}\label{nondimBernoulli}
\frac{q^2}{2}+\int^{\rho}_1\frac{p'(\rho)}{\rho}d\rho
=\frac{1}{2}.
\end{equation}
Since $p'(\rho)/\rho>0$ for $\rho>0$, so (\ref{nondimBernoulli})
yields a representation of the density
\begin{equation}\label{densityspeed}
\rho=g(q^2),
\end{equation}
moreover, $g$ is a decreasing function. For example, for
polytropic gases, after the nondimensionalization,
$p=\rho^{\gamma}/\gamma$, and (\ref{densityspeed}) is nothing but
\begin{equation}
\rho=g(q^2)=\left(\frac{\gamma+1-(\gamma-1)q^2}{2}\right)^{\frac{1}{\gamma-1}}.
\end{equation}

Furthermore, $\rho$ is a two-valued function of $(\rho q)^2$.
Subsonic flows correspond to the branch where $\rho>1$ if $(\rho
q)^2\in[0,1)$. Set
\begin{equation}\label{densitymomentum}
\rho=H((\rho q)^2)
\end{equation}
such that $\rho>1$ if $(\rho q)^2\in[0,1)$, therefore, $H$ is a
positive decreasing function defined on $[0,1]$, twice
differentiable on $[0,1)$, and satisfies $H(1)=1$. Moreover, it
follows from  (\ref{densityspeed}) and (\ref{densitymomentum})
that $(\rho q)^2$ is given in terms of $q^2$ as
\begin{equation}\label{momentumspeed}
(\rho q)^2=G(q^2).
\end{equation}
Thus
\begin{equation*}
g(q^2)=H(G(q^2)).
\end{equation*}

Suppose that the wall of nozzle is impermeable so that
\begin{equation}\label{3Dbc}
(u,v,w)\cdot \overrightarrow{n}=0,
\end{equation}
where $\overrightarrow{n}$ is the outward normal of the solid
boundary.

Due to (\ref{3Dirrotational}), one can introduce a velocity
potential $\Phi$ for the flow such that
\begin{equation*}
\Phi_x=u, \,\, \Phi_y=v,\,\, \Phi_z=w.
\end{equation*}
Thus the continuity equation becomes
\begin{equation*}
\div(g(|\nabla\Phi|^2)\nabla \Phi)=0.
\end{equation*}

Assume now that the nozzle is axi-symmetric as given by
\begin{equation}\label{3Dinteriordomain}
D=\{(x,y,z)|0\leq \sqrt{y^2+z^2}< f(x),-\infty<x<\infty\}.
\end{equation}
Consider a smooth flow in the nozzle. Then it follows from
continuity equation and (\ref{3Dbc}) that mass fluxes through each
section which is transversal to the symmetry axis are the same.

Thus, the problem of finding solutions to smooth flows in a 3-D
nozzle reduces to solving the following problem,
\begin{equation}\label{3Dproblem}
\left\{
\begin{array}{ll}
\div(g(|\nabla\Phi|^2)\nabla\Phi)=0,\qquad &\text{in}\,\, D,\\
\frac{\partial \Phi}{\partial n}=0, \qquad &\text{on}\,\,
\partial D,\\
\int_{S}g(|\nabla\Phi|^2)\frac{\partial\Phi}{\partial \vec{l}}
dS=m_0,
\end{array}
\right.
\end{equation}
where $S$ is the surface transversal to the axis, and $\vec{l}$ is
the normal to $S$ which directs to the positive axial direction.

In this paper, it is assumed that there exists $\alpha\in(0,1)$
such that
\begin{equation}\label{assumption1boundary}
\|f'\|_{C^{\alpha}(\mathbb{R})}<\infty,\,\,
\text{and}\,\,\inf_{\mathbb{R}}f=b>0.
\end{equation}

Now the main results of this paper can be stated as follows.
\begin{thm}\label{Thexistence}
Suppose that the nozzle boundary satisfies
(\ref{assumption1boundary}). Then there exists a positive constant
$\bar{m}$ depending only on $f$ such that if $m_0<\bar{m}$, there
exists an axially symmetric uniformly subsonic flow through the
nozzle. More precisely, there exists a smooth solution $\Phi\in
C^{\infty}(D)$ to (\ref{3Dproblem}) such that
\begin{equation}\label{eqbasicsubsonic}
\sup_{\bar{D}}|\nabla\Phi|< 1,
\end{equation}
and
\begin{equation}\label{3Dbyaxiallysy}
u(x,y,z)=\Phi_x=U(x,r),\,\, v(x,y,z)=\Phi_y=V(x,r)\frac{y}{r},\,\,
w(x,y,z)=\Phi_z=V(x,r)\frac{z}{r},
\end{equation}
where $r=\sqrt{y^2+z^2}$ and $U(x,r)$, $V(x,r)$ are smooth in
their arguments, and $V(x,r)$ vanishes on the symmetry axis.
\end{thm}

If the wall of the nozzle tends to be flat at far fields, for
example, rescaling if necessary, one may assume that
\begin{equation}\label{farflat}
f(x)\rightarrow 1\,\, \text{as}\,\, x\rightarrow -\infty,\,\,
f(x)\rightarrow a>0\,\,\text{as}\,\, x\rightarrow \infty.
\end{equation}
Then the following sharper results hold.
\begin{thm}\label{Thcritical}
Suppose that the wall of the nozzle satisfies both
(\ref{assumption1boundary}) and (\ref{farflat}). Then there exists
$\hat{m}>0$ such that if $0\leq m_0< \hat{m}$, there exists a
unique axially symmetric uniformly subsonic flow through the
nozzle with the properties that
\begin{equation}\label{Criticalvalue}
M(m_0)=\sup_{(x,y,z)\in\bar{D}}\left|\nabla \Phi\right|<1,
\end{equation}
and
\begin{eqnarray*}
&&\left|(U,V)-(\{G^{-1}(\frac{m_0^2}{\pi^2})\}^{1/2},0)\right|\rightarrow
0 \,\,\text{as}\,\,
x\rightarrow -\infty,\\
&&\left|(U,V)-(\{G^{-1}(\frac{m_0^2}{\pi^2
a^4})\}^{1/2},0)\right|\rightarrow 0 \,\,\text{as}\,\,
x\rightarrow \infty,
\end{eqnarray*}
uniformly in $r$, where $G$ is defined by (\ref{momentumspeed});
moreover, $M(m_0)$ ranges over $[0,1)$ as $m_0$ varies in
$[0,\hat{m})$. Furthermore, if $0<m_0<\hat{m}$, the axial velocity
is always positive in $\bar{D}$, i.e.
\begin{equation}\label{positivehorizontal}
u>0,
\end{equation}
and, the flow angle, $\omega=\arctan\frac{V}{U}$, satisfies
\begin{equation}\label{angleestimate}
\underline{\omega}\leq \omega\leq \bar{\omega},
\end{equation}
where
\begin{equation}\label{anglebound}
\underline{\omega}=\min\{\inf_x \arctan f'(x),0\},\,\,
\bar{\omega}=\max\{\sup_x \arctan f'(x),0\}.
\end{equation}
Moreover, for any given $\underline{m}\in (0,\hat{m})$, there
exist a positive constant $\delta=\delta(\underline{m})>0$, such
that if $m\in[\underline{m},\hat{m})$, then
\begin{equation}\label{Tminimumspeedestimate}
q(m)=\inf_{\Omega}|\nabla\Phi|\geq \delta.
\end{equation}
\end{thm}

We now study the limiting behavior of these subsonic flows in the
nozzle when then the cross-section mass fluxes $m_0$ approaches
the critical value. In fact, as $m_0\uparrow \hat{m}$, the
corresponding flow fields tend a limit which yields a
subsonic-sonic flow in the nozzle.

\begin{thm}\label{Thsonic}
Let $\{m_{0,n}\}$ be any sequence such that $m_{0,n}\rightarrow
\hat{m}$ as $n\rightarrow +\infty$. Denote by $(U_n,V_n)$ the
global uniformly subsonic flow corresponding to $m_{0,n}$ as
guaranteed by Theorem 2. Then there exists a subsequence, still
labelled by $\{(U_n,V_n)\}$ associated with $\{m_{0,n}\}$ such
that
\begin{eqnarray}
&U_n\rightarrow U,\quad V_n\rightarrow V,\label{velocityconvegence}\\
&g(q_n^2)U_n\rightarrow g(q^2)U,\quad g(q_n^2)V_n\rightarrow
g(q^2)V,\label{momentumconvergence}
\end{eqnarray}
where $q_n^2=U_n^2+V_n^2$, $q^2=U^2+V^2$, and $g(q^2)$ is the
function defined by (\ref{densityspeed}) through Bernoulli's law,
all the above convergence are almost convergence. Moreover, this
limit yields a three dimensional flow with density
$\rho(x,y,z)=g(q^2)(x,r)$ and velocity
\begin{equation*}
u(x,y,z)=U(x,r),
\,\,v(x,y,z)=V(x,r)\frac{y}{r},\,\,w(x,y,z)=V(x,r)\frac{z}{r},
\end{equation*}
where $r=\sqrt{y^2+z^2}$, which satisfies
\begin{equation*}
u_y-v_x=0,\,\,v_z-w_y=0,\,\,w_x-u_z=0\,\,\text{in}\,\,D
\end{equation*}
in the sense of distribution, moreover, for any $\eta\in
C_c^{\infty}(\bar{D})$
\begin{equation*}
\iiint_{D} (\rho u, \rho v, \rho w)\cdot\nabla\eta dx dydz=0.
\end{equation*}
This implies that $(u,v,w)$ satisfies boundary condition
(\ref{3Dbc}) as the normal trace of the divergence field $(\rho u,
\rho v, \rho w)$ on the boundary.
\end{thm}

Before we prove the theorems, there are a few remarks in order.

\begin{rmk}
{\rm In contrast to two dimensional plane flows, three dimensional
flows are much more complicated. Indeed, some of the key arguments
in \cite{XX} can not be applied to three dimensional case
directly. Even for irrotational steady axially symmetric subsonic
flows, there are some difficulties near the symmetry axis, see
(\ref{streameq}). Therefore, it seems difficult to show the
existence of subsonic flows by fixed point argument as in plane
flows in \cite{XX}. Fortunately, for axisymmetric flows, equation
(\ref{streameq}) has a variational structure, which is one of the
key points to show the existence of subsonic solutions.}
\end{rmk}

\begin{rmk}
{\rm  It should be noted that one cannot adapt the analysis of
\cite{ACF} directly to study the properties of the subsonic flow
in Theorem \ref{Thcritical} since for jet flow, the pressure is
prescribed on the jet surface, so the flow speed is known by
Bernoulli's law, thus it is easier to see whether the flow is
subsonic and whether it approaches to uniform flows at far
fields.}
\end{rmk}

\begin{rmk}
{\rm In all the theorems in this chapter, we require only
$C^{1,\alpha}$ smoothness of $f$. Similar to the proofs given in
this paper, one can prove all results in \cite{XX} under the
condition that nozzle boundaries are $C^{1,\alpha}$ instead of
$C^{2,\alpha}_{loc}$. Furthermore, it is only required that $f$
itself tends to constants at far fields instead of its higher
derivatives, which improves the results in \cite{XX}.}
\end{rmk}

\begin{rmk}
{\rm  Theorem \ref{Thcritical} provides the existence of flows
studied by Gilbarg in \cite{G}. Moreover, applying the comparison
principle obtained by Gilbarg in \cite{G}, the maximum speed of
flows obtained in Theorem \ref{Thcritical} is monotone increasing
with respect to incoming mass flux.}
\end{rmk}

\begin{rmk}
{\rm There are some fragmentary descriptions of some phenomena on
the axially symmetric subsonic flows past a body, for the
reference, please refer to \cite{Bers1}, \cite{GS},
\cite{GSerrin}. For applications of the theory of compensated
compactness to two dimensional transonic and subsonic-sonic flows,
please see \cite{Morawetz}, \cite{CDSW}, \cite{XX}.}
\end{rmk}

The rest of the paper is organized as follows: in Section 2, we
derive the governing equation and boundary conditions for axially
symmetric irrotational flows. In Section 3, we adapt the
variational method used in \cite{ACF} to prove Theorem 1.
Subsequently, in Section 4, we prove that subsonic flows will
approach uniform flows at far fields when the nozzle boundaries
tend to be flat at far fields, which will yield the existence of
the critical value for incoming mass fluxes. In Section 5,
positivity of axial velocity and uniform estimates for flow angles
for axially symmetric flows are proved. In last section, Section
6, we use a compensated compactness framework to show the
existence of weak subsonic-sonic flows.

\section{Axially Symmetric Flows}
In this section, we will derive the governing equations and
boundary conditions for axially symmetric irrotational flows in
cylindrical coordinates and in terms of stream functions.

In the cylindrical coordinates $(x, r,\theta)$, let the fluid
density and velocity be $\rho(x,r,\theta)$ and $(U(x,r,\theta),
V(x,r,\theta), W(x,r,\theta))$, where $U$, $V$, and $W$ are axial
velocity, radial velocity and swirl velocity respectively. Then
$(x,y,z)$, $\rho$, and $(u,v,w)$ satisfy
\begin{eqnarray*}
&& x=x, \,\, y=r\cos\theta,\,\, z=r\sin\theta;\\
&& \rho(x,y,z)=\rho(x,r,\theta), \,\, u(x,y,z)=U(x,r,\theta);\\
&& v(x,y,z)=V(x,r,\theta)\cos\theta+W(x,r,\theta)(-\sin\theta),\\
&&w(x,y,z)=V(x,r,\theta)\sin\theta+W(x,r,\theta)\cos\theta.
\end{eqnarray*}
It should be noted that for axi-symmetric flows, one has
\begin{equation*}
U(x,r,\theta)=U(x,r),\,\,V(x,r,\theta)=V(x,r),\,\,W(x,r,\theta)=W(x,r).
\end{equation*}
Since the flow is also assumed to be irrotational, one has
\begin{equation*}
v_z-w_y=-\frac{(rW)_r}{r}=0,
\end{equation*}
this implies that
\begin{equation*}
W=\frac{c(x)}{r}.
\end{equation*}
Thus $W\equiv0$ since $W$ is bounded near $r=0$. Therefore, for
axially symmetric irrotational flows, one has
\begin{equation}\label{Cycomponent}
u=U(x,r),\,\, v=V(x,r)\frac{y}{r},\,\, w=V(x,r)\frac{z}{r},\,\,
\text{and}\,\, \rho=\rho(x,r),
\end{equation}
where $r=\sqrt{y^2+z^2}$. Then the continuity equation reduces to
\begin{equation}\label{Cyconeq}
(r\rho U)_x+(r\rho V)_r=0.
\end{equation}
Moreover, the irrotational condition (\ref{3Dirrotational})
changes to
\begin{equation}\label{Cyirrotational}
U_r-V_x=0.
\end{equation}
Bernoulli's law (\ref{nondimBernoulli}) is still of the same form
with $q=\sqrt{U^2+V^2}$.

Due to (\ref{Cyconeq}), one can introduce a stream function
$\psi=\psi(x,r)$ such that
\begin{equation}\label{streamfunction}
\psi_r=r\rho U,\,\, \psi_x=-r\rho V.
\end{equation}
Then Bernoulli's law (\ref{nondimBernoulli}) becomes to
\begin{equation}\label{streamBernoulli}
\frac{1}{2\rho^2}\left|\frac{\nabla
\psi}{r}\right|^2+\int^{\rho}_1\frac{p'(\rho)}{\rho}d\rho=\frac{1}{2}.
\end{equation}
Therefore, it follows from (\ref{densitymomentum}) that $\rho$ can
be represented as
\begin{equation}\label{density}
\rho=H\left(\left|\frac{\nabla\psi}{r}\right|^2\right),
\end{equation}
so the irrotationality (\ref{Cyirrotational}) changes to
\begin{equation}\label{streameq}
\div\left(\left(H\left(\left|\frac{\nabla\psi}{r}\right|^2\right)\right)^{-1}
\frac{\nabla\psi}{r}\right)=0.
\end{equation}

The no-flow boundary condition (\ref{3Dbc}) on the nozzle wall
becomes
\begin{equation}\label{Cybc}
(U,V)\cdot \overrightarrow{N}=0,
\end{equation}
where $\overrightarrow{N}$ is the normal of the curve $r=f(x)$. It
follows from (\ref{Cybc}) that $\psi$ is a constant in each
connected component of the solid boundaries.

Note that for smooth axisymmetric flows in the nozzle, it follows
from (\ref{streamfunction}) that $\psi$ is a constant on the
symmetry axis. Thus $r=0$ is a streamline.

Since the flow is axially symmetric, one may consider only
symmetric part of the domain. Let
\begin{equation}\label{domain}
\Omega=\left\{(x,r)\big|0<r<f(x),-\infty<x<\infty\right\}
\end{equation}
with boundaries
\begin{equation}\label{boundary}
T_1=\left\{(x,r)\big|r=0,\,-\infty<x<\infty\right\},\qquad
T_2=\{(x,r)|r=f(x),\,-\infty<x<\infty\}.
\end{equation}
For convenience, we denote by $D_0$ the three dimensional domain
induced by $\Omega$,
\begin{equation}\label{3Ddomainwithboundary}
D_0=\{(x,y,z)|0<\sqrt{y^2+z^2}< f(x),-\infty<x<\infty\}.
\end{equation}

Then, to study the 3-dimensional problem, (\ref{3Dproblem}), for
axisymmetric flows, one may first study the following
2-dimensional problem
\begin{equation}\label{problem}
\left\{
\begin{array}{ll}
\div\left(\left(H\left(\left|\frac{\nabla\psi}{r}\right|^2\right)\right)^{-1}
\frac{\nabla\psi}{r}\right)=0,
\quad &\text{in}\,\,\Omega, \\
\psi=0,\qquad &\text{on}\,\, T_1,\\
\psi=m=\frac{m_0}{2\pi},\qquad &\text{on}\,\, T_2.
\end{array}
\right.
\end{equation}

\section{Subsonic Flows Associated with Small Incoming Mass Flux}

This section is mainly devoted to the proof of Theorem 1. Our
approach is motivated strongly by the important work \cite{ACF} by
Alt, Caffarelli and Friedman. The proof can be divided into 10
steps.

Step 1. Subsonic truncation and shielding singularity. By direct
calculations, it is easy to find that the derivative of function
$H(s)$ goes to negative infinity as $s\rightarrow 1$. To control
the ellipticity and avoid singularity of $H'$, one may truncate
$H$ as follows
\begin{equation}\label{subtruncation}
\tilde{H}(s)=\left\{
\begin{array}{ll}
H(s), \quad &\text{if}\,\, 0\leq s<\tilde{m}^2,\\
H\left((\frac{\tilde{m}+1}{2})^2\right), \quad &\text{if}\,\,
s\geq (\frac{\tilde{m}+1}{2})^2,
\end{array}
\right.
\end{equation}
where $\tilde{m}<1$ is a constant to be determined,  and
$\tilde{H}$ is a smooth decreasing function. Set
\begin{equation}\label{modifyspeedmomentum}
q^2=s/\tilde{H}^2(s).
\end{equation}
Since $\tilde{H}^2(s)-2\tilde{H}\tilde{H}'(s)s>0$,  we can
represent $s$ as a function of $q^2$, $s=\tilde{G}(q^2)$.
Obviously, $\tilde{G}$ is an increasing function. Define
$\rho=\tilde{g}(q^2)$ as
\begin{equation}\label{modifydensityspeed}
\tilde{g}(q^2)=\tilde{H}(\tilde{G}(q^2)).
\end{equation}
Then it is easy to check that
\begin{equation}\label{potentialellipticity}
\Lambda\geq
\tilde{g}+2q^2\frac{d\tilde{g}}{dq^2}=\frac{\tilde{H}(\tilde{G}(q^2))}
{\tilde{H}(\tilde{G}(q^2))
-2\tilde{H}'(\tilde{G}(q^2))\tilde{G}(q^2)}\geq \nu
\end{equation}
for some positive real numbers $\Lambda$ and  $\nu>0$ which depend
on $\tilde{H}$.

To treat the singularity in the coefficients of the equation
(\ref{streameq}) as $r\rightarrow 0$, one may shield the
singularity by first solving the following problem
\begin{equation}\label{hidsingularity}
\left\{
\begin{array}{ll}
\div\left(\left(\tilde{H}\left(\left|\frac{\nabla\psi}{r+\delta}\right|^2
\right)\right)^{-1}
\frac{\nabla\psi}{r+\delta}\right)
=0,\quad &\text{in}\,\,\Omega, \\
\psi=0,\qquad &\text{on}\,\, T_1,\\
\psi=m,\qquad &\text{on}\,\, T_2.
\end{array}
\right.
\end{equation}

Step 2. Variational problem. The problem (\ref{hidsingularity}) is
a boundary value problem for an elliptic equation in a unbounded
domain, therefore, we use a series of Dirichlet problems in
bounded domains to approximate it. Thus consider first the
following problem
\begin{equation}\label{Lhidproblem}
\left\{
\begin{array}{ll}
\div\left(\left(\tilde{H}\left(\left|\frac{\nabla\psi}{r+\delta}\right|^2
\right)\right)^{-1} \frac{\nabla\psi}{r+\delta}\right) =0,
\quad &\text{in}\,\,\Omega_L, \\
\psi=\frac{r^2}{f^2(x)}m,\qquad &\text{on}\,\,
\partial\Omega_L,
\end{array}
\right.
\end{equation}
where $\Omega_{L}=\left\{(x,r)\big|(x,r)\in \Omega,
|x|<L\right\}$. The problem, (\ref{Lhidproblem}), can be solved by
a variational method. The existence of solution to problem
(\ref{Lhidproblem}) is equivalent to find minimizer
$\psi_L^{\delta}\in \mathcal{A}_L=\{\phi|\phi\in
W^{1,2}(\Omega_L), \phi-\frac{r^2}{f(x)}m\in
W^{1,2}_0(\Omega_L)\}$ for the following minimization problem
\begin{equation}\label{Lminimization}
\mathcal{J}_{L}(\psi_L^{\delta})=\inf_{\phi\in
\mathcal{A}_L}\mathcal{J}_{L}(\phi),
\end{equation}
where
\begin{equation}\label{Lfunctional}
\mathcal{J}_L(\phi)=\int_{\Omega_L}F\left(\left|\frac{\nabla
\phi}{r+\delta}\right|^2\right)(r+\delta)dxdr,
\end{equation}
and $F$ is defined by
\begin{equation}\label{primitivedensity}
F(s)=\int_0^s (\tilde{H}(t))^{-1}dt.
\end{equation}
 Since $\tilde{H}$ is a smooth decreasing
function, therefore,
$\mathcal{F}(p_1,p_2,r)=F\left(\left|\frac{(p_1,p_2)}{r+\delta}\right|^2\right)(r+\delta)$
is a convex function of $p=(p_1,p_2)$, using the standard theory
in calculus of variations, for example, Theorem 2 and Theorem 3 in
Section 8.2 in \cite{Evans}, the problem (\ref{Lminimization}) has
a unique solution since the functional $\mathcal{J}_L$ is also
coercive.

Step 3. Estimates for minimizers. For each $L$, there exists a
unique solution $\psi_L^{\delta}$ of the problem
(\ref{Lminimization}). Each minimizer $\psi_L^{\delta}$ is a weak
solution to the Dirichlet problem of the Euler-Lagrange equation,
(\ref{Lhidproblem}). Then by a weak maximum principle for the
problem (\ref{Lhidproblem}), see Theorem 8.1 in \cite{GT}, one
gets
\begin{equation}\label{roughmaximum}
0\leq \psi_L^{\delta}\leq m \,\, \text{in}\,\, \Omega_L.
\end{equation}
Using Caccioppoli's inequality, both in interior and on the
boundary, and Theorem 6.5 and Theorem 6.8 in \cite{Giusti}, one
can obtain
\begin{equation}\label{Caccioppoli}
\|\nabla\psi_k^{\delta}\|_{L^2(\Omega_L)}\leq C(L, \|f'\|_{C^{0}},
\delta, m), \quad \forall k>2L.
\end{equation}
Then H\"{o}ler estimates for the gradient of minimizers to the
functional (\ref{Lfunctional}), and Theorem 8.6 in \cite{Giusti},
imply that there exists $\alpha_1\in (0,\alpha)$ such that
\begin{equation}\label{Holdergradient}
\|\psi_k^{\delta}\|_{C^{1,\alpha_1}(\Omega_L)}\leq C(L, \alpha_1,
\|f'\|_{C^{\alpha}}, \delta, m), \quad \forall k>2L.
\end{equation}
Moreover, the interior Schauder estimate, Theorem 10.18 in
\cite{Giusti}, shows that for any $\Sigma\subset\subset \Omega_L$,
it holds that
\begin{equation}\label{Schauder}
\|\psi_k^{\delta}\|_{C^{2,\alpha_1}(\Sigma)}\leq C(\Sigma, L,
\alpha_1, \delta, m),\quad\forall k>L.
\end{equation}

To recover the singularity later by taking the limit
$\delta\rightarrow 0+$, we need a more precise estimate than
(\ref{roughmaximum}). Set
\begin{equation}\label{supfunction}
\bar{\psi}=\frac{(r+\delta)^2}{b^2}m,
\end{equation}
where $b$ is defined in (\ref{assumption1boundary}). Then it is
easy to check that $\bar{\psi}$ satisfies the equation
\begin{equation}\label{truncateshieldeq}
\div\left(\left(\tilde{H}\left(\left|\frac{\nabla\psi}{r+\delta}\right|^2
\right)\right)^{-1} \frac{\nabla\psi}{r+\delta}\right)=0.
\end{equation}
Because $\Omega_L$ satisfies a uniform exterior cone condition,
$\psi_L^{\delta}\in C^{0}(\overline{\Omega_L})$ by Theorem 8.29 in
\cite{GT}. Moreover, by (\ref{Schauder}), $\psi_L^{\delta}\in
C^{2,\alpha_1}(\Omega_L)$. Therefore, both $\psi_L^{\delta}$ and
$\bar{\psi}$ satisfy the equation in (\ref{Lhidproblem}) on
$\Omega_L$ in the classical sense. Obviously, $\bar{\psi}\geq
\psi_L^{\delta}$ on $\partial \Omega_L$. Thus, it follows from a
comparison principle, Theorem 10.1 in \cite{GT}, that
\begin{equation}\label{refineupperbound}
\psi_L^{\delta}\leq \bar{\psi} \quad \text{in}\,\, \Omega_L.
\end{equation}

Step 4. Existence of solutions to (\ref{hidsingularity}). By a
diagonal process and Arzela-Ascoli lemma, it follows from
(\ref{Holdergradient}) that there exists a sequence $\{n_k\}$ such
that
\begin{equation*}
\psi_{n_k}^{\delta}\chi_{\Omega_{n_k}}\rightarrow
\psi^{\delta}\qquad \text{in}\,\, C^{1,\mu}(\Omega_L)\,\,
\text{for}\,\, \forall L>0
\end{equation*}
with $0<\mu<\alpha_1$. Therefore, $\psi^{\delta}$ is a weak
solution to the problem (\ref{hidsingularity}). Then it follows
from (\ref{Schauder}) that $\psi^{\delta}\in C^{2,\mu}(\Omega_L)$,
$\forall L>0$.

Step 5. Recover singular coefficients. Due to
(\ref{refineupperbound}), we have
\begin{equation*}
\psi^{\delta}(x,r)\leq \frac{m}{b^2}(r+\delta)^2.
\end{equation*}
Therefore, $\psi^{\delta}\rightarrow 0$ on $T_1$. Moreover, for
$\forall \varepsilon>0$, on each set
$\Omega_{L,\varepsilon}=\left\{(x,r)\big||x|<L, \varepsilon
<r<f(x)\right\}$, it follows from Caccioppoli's inequality and
H\"{o}lder gradient estimate in a similar way as for
(\ref{Caccioppoli}) and (\ref{Holdergradient}), that
$\psi^{\delta}$ satisfies the following estimate
\begin{equation}
\|\psi^{\delta}\|_{C^{1,\alpha_1}(\overline{\Omega_{L,
\varepsilon}})}\leq C(L, \varepsilon, \|f'\|_{C^{\alpha}}, m).
\end{equation}
Due to a diagonal process and Arzela-Ascoli Lemma again, there
exists a subsequence $\{\delta_k\}$ such that
\begin{equation*}
\psi^{\delta_k}\rightarrow \psi\,\, \text{in}\,\,
C^{1,\mu}(\Omega_{L,\varepsilon})\,\, \text{for each}\,\, L>0,\,\,
\varepsilon >0.
\end{equation*}
In particular,
\begin{equation}\label{pointwiselimit}
\psi^{\delta_k}\rightarrow \psi\,\, \text{pointwise in}\,\,
\Omega.
\end{equation}
Moreover, $\psi\in C^{1,\mu}$ solves the problem
\begin{equation}\label{truncatepb}
\left\{
\begin{array}{ll}
\div\left(\left(\tilde{H}\left(\left|\frac{\nabla\psi}{r}\right|^2
\right)\right)^{-1} \frac{\nabla\psi}{r}\right)=0,
\quad &\text{in}\,\,\Omega, \\
\psi=\frac{r^2}{f^2(x)}m,\qquad &\text{on}\,\,
\partial\Omega,
\end{array}
\right.
\end{equation}
weakly and satisfies
\begin{equation}\label{upperbound}
0\leq \psi\leq \frac{m}{b^2}r^2.
\end{equation}
It follows from the standard bootstrap arguments that $\psi\in
C^{2,\mu}(\Omega)$.

Step 6. Subsonic estimate near the symmetry axis. In this step,
our aim is to show that
\begin{equation}\label{axissubestimate}
\left|\frac{\nabla \psi}{r}(x,r)\right|\leq Cm
\end{equation}
for $0<r<\frac{b}{2}$. To do this, we note an important
observation due to \cite{ACF}, that if $r_0<b/2$, then
\begin{equation}\label{axiscaling}
\psi_0(x,r)=\frac{1}{t^2}\psi(x_0+tx,r_0+tr), \quad
t=\frac{r_0}{2},
\end{equation}
satisfies
\begin{equation*}
\div\left(\left(\tilde{H}\left(\left|\frac{\nabla\psi_0}{2+r}\right|^2
\right)\right)^{-1} \frac{\nabla\psi_0}{2+r}\right)=0 \quad
\text{in}\,\,B_1\left((0,0)\right).
\end{equation*}
It follows from (\ref{upperbound}) that
\begin{equation*}
0\leq \psi_0\leq Cm \quad\text{in}\,\, B_1\left((0,0)\right).
\end{equation*}
Therefore, by Moser's iteration, Theorem 8.18 in \cite{GT}, one
can get
\begin{equation*}
|\nabla \psi_0|\leq Cm \quad \text{in}\,\,
B_{\frac{1}{2}}\left((0,0)\right).
\end{equation*}
In particular,
\begin{equation}\label{pointaxissubestimate}
\left|\frac{\nabla \psi}{r}(x_0,
r_0)\right|=|\nabla\psi_0(0,0)|\leq Cm.
\end{equation}

Step 7. Subsonic estimate away from the symmetry axis. In this
step, we derive the estimate (\ref{axissubestimate}) for $r>b/4$.
For any given $(x_0,r_0)\in \Omega_{\infty,
b/4}=\{(x,r)|(x,r)\in\Omega, r>b/4\}$, noting that all the
coefficients in the equation (\ref{truncatepb}) are bounded on
$B_{b/8}((x_0,r_0))\bigcap \Omega_{\infty,b/4}$, moreover,
\begin{equation*}
0\leq \psi\leq m
\end{equation*}
due to (\ref{roughmaximum}) and (\ref{pointwiselimit}), one can
derive by Moser's iteration that
\begin{equation*}
|\nabla \psi(x_0,r_0)|\leq Cm.
\end{equation*}
Therefore,
\begin{equation}\label{pointaaxissubestimate}
\left|\frac{\nabla \psi}{r}(x_0, r_0)\right|\leq Cm,
\end{equation}
since $r_0>b/4$.

Step 8. Uniform H\"{o}lder continuity of velocity field near the
symmetry axis. It follows from Step 6 and Step 7 that the flow is
subsonic except on the axis when the incoming mass flux is
sufficiently small. To show that the flow is subsonic globally, we
first need to show that the velocity field is well-defined along
the symmetry axis. In fact, we have the following stronger
results.
\begin{lem}\label{lemaxisestimate}
Let $\psi$ be a solution to problem (\ref{truncatepb}) satisfying
\begin{equation}\label{lemboundspeedassumption}
\left|\frac{\nabla\psi}{r}\right|\leq C\qquad \text{in}\,\,\Omega.
\end{equation}
Then $\frac{\nabla\psi}{r}$ is uniformly H\"{o}lder continuous up
to the symmetry axis, moreover,
\begin{equation}\label{velocityonaxis}
\lim_{(x,r)\rightarrow (x_0,0)}\frac{\psi_x}{r}(x,r)=0
\,\,\text{for any}\,\, x_0\in(-\infty,\infty).
\end{equation}
More precisely, there exists $\beta\in(0,1)$ such that
\begin{equation}\label{Holderaxis}
\left[\frac{\nabla
\psi}{r}\right]_{C^{\beta}((l_1,l_2)\times(0,h_0))}\leq
C(|l_2-l_1|)
\end{equation}
holds for $0<h_0\leq \frac{b}{4}$ and any real numbers $l_1<l_2$.
\end{lem}
\begin{pf}
Since $\psi\in C^{2,\mu}(\Omega)\bigcap C^{1,\mu}(\Omega\bigcup
T_2)$  satisfies the equation in (\ref{truncatepb}) and admits the
bound (\ref{lemboundspeedassumption}), therefore, the axially
symmetric potential
\begin{equation}\label{axialsypotential}
\varphi(x,r)=\int_{(0,f(0))}^{(x,r)}\left(
\tilde{H}\left(\left|\frac{\nabla\psi}{r}\right|^2\right)\right)^{-1}
\frac{\psi_r}{r}dx
-\left(\tilde{H}\left(\left|\frac{\nabla\psi}{r}\right|^2\right)\right)^{-1}
\frac{\psi_x}{r}dr
\end{equation}
is well-defined and path independent except on the symmetry axis
$\{r=0\}$. Moreover,
\begin{equation}\label{sypotentialderivative}
\varphi_x=\left(\tilde{H}\left(\left|\frac{\nabla\psi}{r}\right|^2\right)\right)^{-1}
\frac{\psi_r}{r},\,\,
\varphi_r=-\left(\tilde{H}\left(\left|\frac{\nabla\psi}{r}\right|^2\right)\right)^{-1}\frac{\psi_x}{r}.
\end{equation}
Therefore, by (\ref{lemboundspeedassumption}),
\begin{equation}\label{boundedsypotentialderivative}
|\nabla \varphi|\leq C \,\,\text{in}\,\, \Omega=\{0<r<f(x)\},
\end{equation}
which implies that $\varphi$ can be extended to $\bar{\Omega}$ as
\begin{equation*}
\varphi(x_0,0)=\lim_{(x,r)\in\Omega,\\ (x,r)\rightarrow
(x_0,0)}\varphi(x,r),\qquad \forall x_0\in(-\infty,\infty).
\end{equation*}
Note that the axially symmetric potential $\varphi$ induces a 3-D
potential function
\begin{equation}\label{def3Dpotential}
\Phi(x,y,z)=\varphi(x,\sqrt{y^2+z^2})
\end{equation}
which is defined on the three dimensional domain $\bar{D}$. Then,
\begin{equation}\label{3Dpotentialderivative}
\Phi_x=\varphi_x,\,\,
\Phi_y=\varphi_r\frac{y}{r},\,\,\Phi_z=\varphi_r\frac{z}{r}\qquad
\text{in}\,\,D_0,
\end{equation}
and that $\Phi\in W_{loc}^{1,\infty}(D)\bigcap C^{2,\mu}(D_0)$.
Moreover, it follows from (\ref{sypotentialderivative}),
(\ref{3Dpotentialderivative}), and (\ref{modifydensityspeed}) that
for any $\varepsilon>0$, $\Phi$ solves the equation
\begin{equation}\label{modify3Dpotentialeq}
\div (\tilde{g}(|\nabla\Phi|^2)\nabla \Phi)=0
\end{equation}
in the three dimensional domain
$D_{\varepsilon}=\{(x,y,z)|-\infty<x<\infty, \varepsilon
<\sqrt{y^2+z^2}<f(x)\}$. Therefore, any $\eta\in
C_{0}^{\infty}(D)$, one has
\begin{eqnarray*}
&&\iiint_{D}\tilde{g}(|\nabla \Phi|^2)\nabla\Phi\cdot\nabla\eta dx
dydz\\
&=&\iiint_{D_{\varepsilon}}\tilde{g}(|\nabla
\Phi|^2)\nabla\Phi\cdot\nabla\eta dx dy dz +\iiint_{D \setminus
D_{\varepsilon}}\tilde{g}(|\nabla
\Phi|^2)\nabla\Phi\cdot\nabla\eta dx dy
dz\\
&=& -\iint_{\partial
D_{\varepsilon}}\tilde{g}(|\nabla\Phi|^2)\frac{\partial\Phi}{\partial
n}\eta dS + \iiint_{D \setminus D_{\varepsilon}}g(|\nabla
\Phi|^2)\nabla\Phi\cdot\nabla\eta dx dy
dz\\
&=&-\int_{-\infty}^{\infty}\int_0^{2\pi}\tilde{g}(|\nabla\Phi|^2)
\frac{\partial}{\partial r}\Phi(x,\varepsilon,\theta)\eta
\varepsilon dxd\theta +
\int_{-\infty}^{\infty}\int_0^{\varepsilon}\int_0^{2\pi} g(|\nabla
\Phi|^2)\nabla\Phi\cdot\nabla\eta r dx drd\theta\\
 &\rightarrow & 0
\end{eqnarray*}
as $\varepsilon \rightarrow 0$ since $\nabla\Phi$ is bounded.
Therefore,
\begin{equation}\label{weakform3Dpotential}
\iiint_{D}\tilde{g}(|\nabla \Phi|^2)\nabla\Phi\cdot\nabla\eta dx
dydz=0\qquad \forall \eta\in C_{0}^{\infty}(D).
\end{equation}
Thus, $\Phi$ is a weak solution of equation
(\ref{modify3Dpotentialeq}) in $D$.
Since(\ref{modify3Dpotentialeq}) is elliptic due to
(\ref{potentialellipticity}), thus the standard elliptic
regularity theory, \cite{GT}, shows
\begin{equation*}
\Phi\in C^{\infty}(D).
\end{equation*}
Moreover, for $k=1,2,3$, $\partial_k\Phi$ satisfies the equation
\begin{equation}\label{mpotentialderivativeeq}
\partial_i\left((\tilde{g}(|\nabla\Phi|^2)\delta_{ij}
+2\tilde{g}'(|\nabla\Phi|^2)\partial_i\Phi\partial_j\Phi)
\partial_j(\partial_k\Phi)\right)=0,
\end{equation}
which is uniformly elliptic due to (\ref{potentialellipticity}).
Thus, by Nash-Moser iteration, there exists a $\beta_1\in(0,1)$
such that for suitably small positive constant $h$,
\begin{equation}\label{3Dpotentialderivativeestimate}
[\partial_k\Phi]_{C^{\beta_1}(B_h((x,0,0)))}\leq C\|\nabla
\Phi\|_{L^{\infty}(B_{2h}((x,0,0)))}.
\end{equation}
Therefore,
\begin{equation*}
|\partial_y\Phi(x,y,z)-\partial_y\Phi(x,0,0)|\leq
C(y^2+z^2)^{\beta_1/2}\,\,\text{for}\,\, r=(y^2+z^2)^{1/2}\leq h.
\end{equation*}
It follows from (\ref{3Dpotentialderivative}) that
\begin{equation*}
|\varphi_r(x,r)\cos\theta-\Phi_y(x,0,0)|\leq
Cr^{\beta_1},\,\,\text{for all}\,\,\theta\in[0,2\pi).
\end{equation*}
Thus,
\begin{equation*}
\Phi_y(x,0,0)=0.
\end{equation*}
Similarly, $\Phi_z(x,0,0)=0$. Therefore
\begin{equation*}
|\varphi_r(x,r)|\leq Cr^{\beta_1}.
\end{equation*}
Thus,
\begin{equation}\label{sypotetialzeroradial}
\lim_{(x,r)\rightarrow (x_0,0)}\varphi_r(x,r)=0.
\end{equation}
Furthermore, (\ref{3Dpotentialderivativeestimate}) yields
\begin{equation}\label{sypotentialholder}
[\nabla \varphi]_{C^{\beta_1}((l_1,l_2)\times(0,h))}\leq
C(|l_2-l_1|).
\end{equation}
So the desired estimates (\ref{velocityonaxis}) and
(\ref{Holderaxis}) follow.

This finishes the proof of the Lemma.
\end{pf}

Step 9. Removal of cutoff. Combining (\ref{pointaxissubestimate}),
 (\ref{pointaaxissubestimate}) and (\ref{Holderaxis}) yields
\begin{equation*}
\left|\frac{\nabla \psi}{r}\right|\leq Cm,
\,\,\text{in}\,\,\bar{\Omega}.
\end{equation*}
If $m$ is sufficiently small, then $Cm<\tilde{m}$, therefore
\begin{equation*}
\left|\frac{\nabla \psi}{r}\right|\leq \tilde{m}.
\end{equation*}
Consequently, $\psi$ solves the problem (\ref{problem}), and
moreover, which is uniformly subsonic.

Step 10. Existence of 3-D subsonic flow. It follows from the proof
of Lemma \ref{lemaxisestimate} and step 1-9 that there exists a
three dimensional subsonic solution to problem (\ref{3Dproblem})
which satisfies (\ref{eqbasicsubsonic}) and (\ref{3Dbyaxiallysy}).

\section{Existence of The Critical Incoming Mass Flux}

In this section, it will be shown that there exists a critical
value $\hat{m}$ such that the flow is always subsonic when the
three dimensional mass flux $m_0$ is less than $\hat{m}$. To
achieve this goal, we first show that the flow approximates to
uniform flows at far fields.

Let
\begin{equation}\label{secondmodifysvolume}
\hat{H}(s)=\left\{
\begin{array}{ll}
H(s)\qquad &\text{if},\,\, s<s_0^2,\\
H((\frac{s_0+1}{2})^2),\qquad &\text{if} \,\,
s>(\frac{s_0+1}{2})^2
\end{array}
\right.
\end{equation}
be a smooth decreasing function, where $s_0\in (0,1)$. It follows
from the proof of Theorem \ref{Thexistence} that  there exists a
solution $\psi$ to the problem
\begin{equation}\label{secondmodifyproblem}
\left\{
\begin{array}{ll}
\div\left(\left(\hat{H}\left(\left|\frac{\nabla\psi}{r}\right|^2
\right)\right)^{-1}
\frac{\nabla\psi}{r}\right)=0,
\quad &\text{in}\,\,\Omega, \\
\psi=0,\qquad &\text{on}\,\, T_1,\\
\psi=m,\qquad &\text{on}\,\, T_2,
\end{array}
\right.
\end{equation}
for any $m>0$. Moreover, $\psi$ satisfies
\begin{equation}\label{orderforstream}
\psi\leq Cr^2.
\end{equation}
If the wall of the nozzle tends to be flat at far fields, i.e.,
$f$ satisfies (\ref{farflat}), then solutions to
(\ref{secondmodifyproblem}) approximate to uniform flows at far
fields, as is described in the following lemma.
\begin{lem}\label{lemuniformflow}
Suppose that $f$ satisfies (\ref{assumption1boundary}) and
(\ref{farflat}). Let $\psi$ be a solution to
(\ref{secondmodifyproblem}) and satisfy (\ref{orderforstream}).
Then for any $\varepsilon>0$, there exists a constant $L>0$ such
that
\begin{equation*}
\left|\frac{\nabla\psi}{r}(x,r)-(0,2m)\right|<\varepsilon,\,\,
\text{if}\,\,x<-L,
\end{equation*}
and
\begin{equation*}
\left|\frac{\nabla\psi}{r}(x,r)-(0,\frac{2m}{a^2})\right|<\varepsilon,\,\,\text{if}\,\,
x>L.
\end{equation*}
\end{lem}
\begin{pf}
We start with a special case. Assume that $f(x)=a$ if $x>L_0$. Set
$\psi_k(x,r)=\psi(x+k,r)\chi_{\{(x,r)| x>-k+L_0+1, 0<r<a\}}$. It
follows from the standard H\"{o}lder gradient estimates and Lemma
\ref{lemaxisestimate} that there exists $\alpha_2\in(0,1)$ such
that for any compact set $K\subset (-\infty,\infty)\times [0,a]$,
$\|\psi_k\|_{C^{1,\alpha_2}(K)}\leq C(K)$ for $k$ sufficiently
large, where $C(K)$ does not depend on $k$. Therefore, by
Arzela-Ascoli lemma, there exists a subsequence
$\psi_{k_l}\rightarrow \psi_0$ in $C^{1,\alpha_3}(K)$ with
$\alpha_3<\alpha_2$. Moreover, $\psi_0$ solves the following
boundary value problem
\begin{equation}\label{blowupinfinitypb}
\left\{
\begin{array}{ll}
\mathcal{L}\psi_0=\div
\left((\hat{H}(|\frac{\nabla\psi_0}{r}|^2))^{-1}\frac{\nabla\psi_0}{r}\right)=0,
\qquad &\text{in}\,\,
E_0=\{(x,r)|-\infty<x<\infty, 0<r<a\},\\
\psi_0(x,a)=m,\qquad &\text{if}\,\,-\infty<x<\infty,\\
\psi_0(x,0)=0,\qquad &\text{if}\,\,-\infty<x<\infty.
\end{array}
\right.
\end{equation}
Furthermore, thanks to (\ref{orderforstream}), $\psi_0$ satisfies
the estimate
\begin{equation}\label{blowuporderforstream}
\psi_0\leq Cr^2\qquad \text{in}\,\, E_0.
\end{equation}

In fact, problem (\ref{blowupinfinitypb}) and
(\ref{blowuporderforstream}) has a unique solution
\begin{equation}\label{blowupinfinityfunction}
\psi_0=\frac{mr^2}{a^2},
\end{equation}
which follows from a simple comparison argument in \cite{ACF}.

Since the solution to problem (\ref{blowupinfinitypb}) and
(\ref{blowuporderforstream}) is unique, therefore, for
$K=[-2,2]\times [h,a]$, then
\begin{equation*}
\left\|\frac{\nabla \psi_k}{r}-\frac{\nabla
\psi_0}{r}\right\|_{C^{\mu}([-2,2]\times [h,a])}\rightarrow
0\qquad \text{for}\,\, \forall h>0,\,\,
\text{for}\,\,\mu<\alpha_3.
\end{equation*}
By the definition of $\psi_k$ and (\ref{blowupinfinityfunction}),
this is equivalent to
\begin{equation*}
\left\|\frac{\nabla
\psi}{r}-(0,2m/a^2)\right\|_{C^{\mu}([k-2,k+2]\times
[h,a])}\rightarrow 0\,\, \text{as}\,\, k\rightarrow\infty \,\,
\text{for}\,\, \forall h>0\,\, \text{for}\,\,\mu<\alpha_3.
\end{equation*}

In the general case that the wall of the nozzle is not flat at far
fields, one can set
$\psi_k(x,r)=\psi(x+k,r)\chi_{\{(x,r)|x>-k+1,0<r<f(x+k)\}}$. Then
it follows from a similar analysis that $\psi_{k_l}\rightarrow
\psi_0$ in $C^{2,\alpha_3}(K)$ for any compact set $K\subset
(-\infty,\infty)\times(0,a)$, here $K$ may not touch the boundary
$r=a$, but $\psi_0$ still satisfies the same boundary value
problem (\ref{blowupinfinitypb}) and estimate
(\ref{blowuporderforstream}). Therefore,
\begin{equation}\label{interfarfieldestimate}
\left\|\frac{\nabla
\psi}{r}-(0,2m/a^2)\right\|_{C^{\mu}([k-2,k+2]\times
[h,a-h])}\rightarrow 0\,\, \text{as}\,\, k\rightarrow \infty\,\,
\text{for}\,\, \forall h>0 \,\, \text{for}\,\,\mu<\alpha_3.
\end{equation}
However, away from the symmetry axis, $\psi$ possesses H\"{o}lder
gradient estimates, consequently, there exists $\alpha_4>0$ such
that
\begin{equation}\label{boundaryfarfieldestimate}
\left[\frac{\nabla
\psi}{r}\right]_{C^{\alpha_4}(\{(x,r)|k-2<x<k+2, r>a-h\})}\leq C.
\end{equation}
Following the same argument in Section 3, one can show that there
exists $\beta_2\in (0,1/2)$ such that
\begin{equation}\label{axisfarfieldestimate}
\left[\frac{\nabla\psi}{r}\right]_{C^{\beta_2}((k-2,k+2)\times[0,h])}\leq
C.
\end{equation}

It now follows from estimates (\ref{interfarfieldestimate}),
(\ref{boundaryfarfieldestimate}) and (\ref{axisfarfieldestimate})
that the flow approximates to uniform flows at far fields. Indeed,
by (\ref{boundaryfarfieldestimate}) and
(\ref{axisfarfieldestimate}),
\begin{equation}\label{farfieldnoninteriorestimate}
\left|\frac{\nabla
\psi}{r}(x,r_1)-\frac{\nabla\psi}{r}(x,r_2)\right|\leq
C(h_0^{\beta_2}+h_0^{\alpha_4})\,\,
\text{if}\,\,r_1,r_2>a-h_0\,\,\text{or}\,\, 0<r_1,r_2<h_0.
\end{equation}
Thus, $\forall \varepsilon>0$, it follows from
(\ref{farfieldnoninteriorestimate}) that there exists $\bar{h}>0$
such that
\begin{equation}\label{uniformaxisestimate}
\left|\frac{\nabla
\psi}{r}(x,r_1)-\frac{\nabla\psi}{r}(x,r_2)\right|\leq
\frac{\varepsilon}{3},\qquad \forall 0\leq r_1, r_2<\bar{h}.
\end{equation}
and
\begin{equation}\label{uniformsolidestimate}
\left|\frac{\nabla
\psi}{r}(x,r_1)-\frac{\nabla\psi}{r}(x,r_2)\right|\leq
\frac{\varepsilon}{3},\qquad \forall r_1, r_2>a-\bar{h}.
\end{equation}
On the other hand, there exists $L>0$ such that
\begin{equation}\label{uniforminteriorestimate}
\left|\frac{\nabla \psi}{r}(x,r)-(0,2m/a^2)\right|\leq
\frac{\varepsilon}{3},\qquad \frac{\bar{h}}{2}<
r<a-\frac{\bar{h}}{2}, x>L.
\end{equation}
Thus, combining (\ref{uniformaxisestimate}),
(\ref{uniformsolidestimate}) and (\ref{uniforminteriorestimate})
yields
\begin{equation*}
\left|\frac{\nabla \psi}{r}(x,r)-(0,\frac{2m}{a^2})\right|\leq
\varepsilon,\qquad \forall x>L.
\end{equation*}
Similarly, we have
\begin{equation*}
\left|\frac{\nabla \psi}{r}(x,r)-(0,2m)\right|\leq
\varepsilon,\qquad \forall x<-L.
\end{equation*}

This implies that the flow approximates to uniform flows at far
fields.
\end{pf}

With the help of Lemma \ref{lemuniformflow}, one can show the
uniqueness of uniformly subsonic flows.
\begin{lem}\label{lemuniqueness}
Suppose that $f$ satisfies (\ref{assumption1boundary}) and
(\ref{farflat}), then uniformly subsonic flows to problem
(\ref{problem}) are unique.
\end{lem}
\begin{pf} The proof is quite similar to the proof in \cite{XX}.
Suppose there are two uniformly subsonic flows $\psi_1$ and
$\psi_2$ which satisfy
\begin{equation*}
\left|\frac{\nabla \psi_1}{r}\right|,\,\, \left|\frac{\nabla
\psi_2}{r}\right|\leq s_0<1.
\end{equation*}
Since $\psi_i$ ($i=1,2$) satisfies the equation
\begin{equation*}
\div\left(\left(\hat{H}\left(\left|\frac{\nabla\psi_i}{r}\right|^2
\right)\right)^{-1}
\frac{\nabla\psi_i}{r}\right)=0
\end{equation*}
where $\hat{H}$ is defined in (\ref{secondmodifysvolume}). It is
easy to check that $\bar{\psi}=\psi_1-\psi_2$ satisfies an
equation of the form
\begin{equation*}
\tilde{\mathcal{L}}\bar{\psi}=a_{ij}(x,r)\partial_{ij}\bar{\psi}
+b_i(x,r)\partial_i\bar{\psi}=0.
\end{equation*}
By Lemma \ref{lemuniformflow}, the flows corresponding to $\psi_1$
and $\psi_2$ approximate to same uniform flows at the far fields,
therefore, for any $\varepsilon >0$, there exists a $L>0$ such
that $|\bar{\psi}(x,r)|<\varepsilon$ if $|x|>L$. Thus by maximum
principle, $|\bar{\psi}|<\varepsilon$, $\forall\varepsilon>0$,
since $\bar{\psi}=0$ on $T_1$ and $T_2$. Since $\varepsilon$ is
arbitrary, so $\bar{\psi}=0$.

This finishes the proof of the Lemma.
\end{pf}

With the help of Lemma \ref{lemuniformflow}, Lemma
\ref{lemuniqueness}, and going back to the original three
dimensional flows, we can show in the same way as in \cite{XX}
that there exists $\hat{m}$ such that as $m_0\rightarrow \hat{m}$,
$M(m_0)\rightarrow 1$. Furthermore, it follows from the comparison
principle by Gilbarg\cite{G} that as $m_0\uparrow \hat{m}$,
$M(m_0)\uparrow 1$.

\section{Properties of Subsonic Flows}

In this section, as same as the case for plane flows, we will
obtain some properties of subsonic axially symmetric flows, which
are useful to show the existence of subsonic-sonic flows.

It follows from the proof of Lemma \ref{lemaxisestimate} that for
axially symmetric subsonic flows, problem (\ref{3Dproblem}) and
(\ref{problem}) are equivalent. Thus, in this section, we will use
two descriptions simultaneously.

First of all, as for the plane flows in \cite{XX}, the axial
velocity is always positive.
\begin{lem}\label{lempositiveaxialvelocity}
Suppose that $f$ satisfies (\ref{assumption1boundary}) and
(\ref{farflat}). Let $\psi$ be a uniformly subsonic solution to
problem (\ref{problem}), then
\begin{equation*}
u>0\qquad \text{in}\,\,\bar{D}.
\end{equation*}
\end{lem}
\begin{pf}
Since the flow is uniformly subsonic, one can assume that
\begin{equation*}
\sup_{\bar{\Omega}}\left|\frac{\nabla \psi}{r}\right|\leq s_0<1.
\end{equation*}
Define $\hat{g}$ as in (\ref{modifydensityspeed}) with the help of
$\hat{H}$ in (\ref{secondmodifysvolume}), and $\Phi$ as in
(\ref{def3Dpotential}). Then $\Phi$ satisfies
\begin{equation*}
\div(\hat{g}(|\nabla \Phi|^2)\nabla\Phi)=a_{ij}(\nabla
\Phi)\partial_{ij}\Phi=0, \qquad \text{in}\,\, D.
\end{equation*}
Set $u=\Phi_x$. Then
\begin{equation*}
a_{ij}(\nabla
\Phi)\partial_{ij}u+D_{p_k}a_{ij}(\nabla\Phi)\partial_{ij}\Phi\partial_k
u=0,
\end{equation*}
which is a uniform elliptic equation of $u$. Since $\psi=m$ on the
solid boundary $r=f(x)$, therefore, $\psi_x+\psi_rf'(x)=0$. On the
other hand, $\psi$ attains its maximum on the solid boundary,
therefore, by Hopf lemma, $\frac{\partial \psi}{\partial
\vec{N}}>0$, where $\vec{N}$ is unit outward normal to the 2-D
domain $\Omega$. Since
\begin{equation*}
\frac{\partial\psi}{\partial
\vec{N}}=\psi_x(-f'(x)/\sqrt{1+(f'(x))^2})+\psi_r/\sqrt{1+(f'(x))^2}=\psi_r\sqrt{1+(f'(x))^2},
\end{equation*}
thus, $\psi_r>0$. By the definition of $\varphi$ and $\Phi$,
$\psi_r>0$ is equivalent to $\Phi_x>0$ on the three dimensional
solid boundary $\{(x,r,z)|f(x)=\sqrt{y^2+z^2}\}$.  Since the flow
approximates to uniform flows at far fields, moreover,
$\Phi_x\rightarrow \{G^{-1}(\frac{4m^2}{a^4})\}^{1/2}$ as
$x\rightarrow \infty$, and $\Phi_x\rightarrow
\{G^{-1}(4m^2)\}^{1/2}$ as $x\rightarrow -\infty$,, therefore, by
the maximum principle
\begin{equation*}
u>0   \qquad \text{in}\,\, \bar{D}.
\end{equation*}
Therefore, the proof of the Lemma is complete.
\end{pf}

Since the axial velocity is positive, then we can define the flow
angle by
\begin{equation}\label{defangle}
\omega=\arctan \frac{V}{U}.
\end{equation}
Moreover, we have the following estimate on the flow angles.
\begin{lem}\label{lemangleestimate}
Suppose that $f$ satisfies (\ref{assumption1boundary}) and
(\ref{farflat}). Let $\psi$ be a uniformly subsonic solution to
the problem (\ref{problem}), then the angle $\omega$ defined by
(\ref{defangle}) satisfies
\begin{equation}\label{eqlemangleestimate}
\underline{\omega}\leq \omega\leq \bar{\omega},
\end{equation}
where
\begin{equation*}
\underline{\omega}=\min\{0,\inf_x\arctan f'(x)\},\,\,
\bar{\omega}=\max\{0,\sup_x\arctan f'(x)\}.
\end{equation*}
\end{lem}
\begin{pf}
The basic idea for the proof of the lemma is the same as that in
\cite{XX}, i.e., using hodograph transformation to obtain an
elliptic equation for the angle, then the estimate
(\ref{eqlemangleestimate}) will be obtained by a comparison
principle for elliptic equations. However, for axially symmetric
flow, this procedure is more involved.

Let us first go back to the equations for axially symmetric flows
(see (\ref{Cyconeq}) and (\ref{Cyirrotational})), which reads
\begin{equation*}
\left\{
\begin{array}{ll}
(r g U)_x+(r g V)_r=0,\\
U_r-V_x=0,
\end{array}
\right.
\end{equation*}
where $U=\frac{\psi_r}{r g}$, $V=-\frac{\psi_x}{r g}$, $g=g(q^2)$,
$q=\sqrt{U^2+V^2}$, and $g$ is defined by (\ref{densityspeed}).
Since $\psi$ satisfies (\ref{problem}), so $\varphi$ as in
(\ref{axialsypotential}) is well-defined, moreover,
\begin{equation*}
J=\frac{\partial(\varphi,\psi)}{\partial(x,r)}=r g(q^2)q^2\geq 0,
\end{equation*}
and which is strictly positive for $r>0$. For $r>0$, the mapping
$(x,r)\mapsto (\varphi,\psi)$ is a local differmorphism. In fact,
the mapping is globally invertible. Indeed, suppose that there are
two points $(x_1,r_1)$ and $(x_2,r_2)$ such that
$\varphi(x_1,r_1)=\varphi(x_2,r_2)$ and
$\psi(x_1,r_1)=\psi(x_2,r_2)$. If $\psi(x_1,r_1)=\psi(x_2,r_2)=0$,
then it is obvious that $x_1=x_2$ and $r_1=r_2=0$ due to maximum
principle and Lemma \ref{lempositiveaxialvelocity}. Let
$\psi(x_1,r_1)=\psi(x_2,r_2)=d>0$, then $(x_1,r_1)$ and
$(x_2,r_2)$ are both on the streamline defined as follows
\begin{equation*}
\left\{
\begin{array}{ll}
\frac{dx}{ds}=U(x,r),\\
\frac{dy}{ds}=V(x,r),\\
x(0)=x_1,r(0)=r_1.
\end{array}
\right.
\end{equation*}
Moreover, this streamline is uniformly away from the symmetry
axis, then it follows from the argument in \cite{XX} that
$x_1=x_2$ and $r_1=r_2$.

Now direct calculations show that
\begin{eqnarray*}
&&r((g(q^2) U)_x+(g(q^2) V)_r)+g(q^2) V=r\left(q g(q^2)
(1-\frac{q^2}{c^2})q_{\varphi}+r
g(q^2)^2q^2\omega_{\psi}\right)+g(q^2)
q\sin\omega,\\
&&U_r-V_x=-q^2\omega_{\varphi}+r g(q^2) qq_{\psi},
\end{eqnarray*}
where $c$ is the sound speed. Therefore
\begin{equation}\label{firstangleequation}
\left(\frac{q}{r
g(q^2)}\omega_{\varphi}\right)_{\varphi}+\left(\frac{r g(q^2)
q}{1-\frac{q^2}{c^2}}\omega_{\psi}\right)_{\psi}+
\left(\frac{1}{r(1-\frac{q^2}{c^2})}\sin\omega\right)_{\psi}=0,
\end{equation}
which can be rewritten as
\begin{eqnarray}
&&(\frac{q}{r g(q^2)}\omega_{\varphi})_{\varphi}+(\frac{r g(q^2)
q}{1-\frac{q^2}{c^2}}\omega_{\psi})_{\psi}+\frac{1}{r}\frac{d}{dq}(1-q^2/c^2)
\frac{\sin \omega}{r g(q^2)} \omega_{\varphi}\nonumber\\
+&&\frac{\cos\omega}{r(1-\frac{q^2}{c^2})}\omega_{\psi}-\frac{\sin\omega}{r^2(1-\frac{q^2}{c^2})}
r_{\psi}=0.\label{angleequation}
\end{eqnarray}
It follows from definitions of $\varphi$ and $\psi$ that
\begin{equation}\label{radialhododerivative}
r_{\psi}=\frac{\varphi_x}{\varphi_x\psi_r-\psi_x\varphi_r}=\frac{U}{r
g(q^2) q^2}.
\end{equation}
Substituting (\ref{radialhododerivative}) into
(\ref{firstangleequation}) yields
\begin{eqnarray}
&&(\frac{q}{r g(q^2)}\omega_{\varphi})_{\varphi}+(\frac{r g(q^2)
q}{1-\frac{q^2}{c^2}}\omega_{\psi})_{\psi}+\frac{1}{r}\frac{d}{dq}(1-q^2/c^2)
\frac{\sin \omega}{r g(q^2)} \omega_{\varphi}\nonumber\\
+&&\frac{\cos\omega}{r(1-\frac{q^2}{c^2})}\omega_{\psi}-\frac{U}{r^3
g(q^2)q^2(1-\frac{q^2}{c^2})}\sin\omega
=0.\label{secondangleequation}
\end{eqnarray}
Since the flow is subsonic, $1-\frac{q^2}{c^2}>0$, therefore,
equation (\ref{secondangleequation}) is an elliptic equation. Note
that, by Lemma \ref{lempositiveaxialvelocity}, $\omega\in
(-\frac{\pi}{2},\frac{\pi}{2})$. Moreover, in the domain
$\Omega^{+}=\{\omega>0\}\bigcap \Omega$, it follows from Lemma
\ref{lempositiveaxialvelocity} that $\omega$ satisfies that
\begin{equation*}
\mathcal{Q}\omega=\left(\frac{q}{r
g(q^2)}\omega_{\varphi}\right)_{\varphi}+\left(\frac{r g(q^2)
q}{1-\frac{q^2}{c^2}}\omega_{\psi}\right)_{\psi}+\frac{1}{r}\frac{d}{dq}(1-q^2/c^2)
\frac{\sin \omega}{r g(q^2)} \omega_{\varphi}
+\frac{\cos\omega}{r(1-\frac{q^2}{c^2})}\omega_{\psi}\geq 0.
\end{equation*}
By the maximum principle, Theorem 3.1 in \cite{GT},
\begin{equation}\label{upperangleestimate}
\sup_{\Omega^{+}}\omega\leq \sup_{\omega\in
\partial\Omega^{+}}\omega=\max\{\sup_{\partial\Omega}\omega, 0\}.
\end{equation}
Similarly, in the domain $\Omega^{-}=\{\omega<0\}\bigcap \Omega$,
$\omega$ satisfies
\begin{equation*}
\mathcal{Q}\omega=\left(\frac{q}{r
g(q^2)}\omega_{\varphi}\right)_{\varphi}+\left(\frac{r g(q^2)
q}{1-\frac{q^2}{c^2}}\omega_{\psi}\right)_{\psi}+\frac{1}{r}\frac{d}{dq}(1-\frac{q^2}{c^2})
\frac{\sin \omega}{r } \omega_{\varphi}
+\frac{\cos\omega}{r(1-\frac{q^2}{c^2})}\omega_{\psi}\leq 0,
\end{equation*}
by the maximum principle,
\begin{equation}\label{lowerangleestimate}
\inf_{\Omega^{-}}\omega \geq \inf_{\omega\in
\partial\Omega^{-}}\omega=\min\{\inf_{\partial\Omega}\omega,0\}.
\end{equation}

Combining (\ref{upperangleestimate}) and
(\ref{lowerangleestimate}) together, we have
\begin{equation*}
\min(\inf_{\partial \Omega} \omega, 0)\leq \omega\leq
\max(\sup_{\partial \Omega} \omega, 0).
\end{equation*}
This finishes the proof of the Lemma.
\end{pf}

At the end of this section, we would like to study the
relationship between flow speed and incoming mass flux.
\begin{lem}\label{Alemmonotonicity}
\begin{enumerate}
\item[(a)] Let $0<m_1\leq m_2<\hat{m}$. Suppose that $\psi_i$ are
uniform subsonic solutions to (\ref{problem}) associated with the
incoming mass flux $m_i$($i=1,2$). Then
\begin{equation}\label{Amonotonicitybdy}
|\nabla \psi_1(x,f(x)|\leq |\nabla \psi_2(x,f(x))|,\quad \forall
x\in\mathbb{R}^1
\end{equation}
\item[(b)] Both the supremum of flow speed and infimum of horizontal velocity
for uniformly subsonic solutions to (\ref{problem}) are monotone
increasing with respect to the incoming mass flux. Moreover, for any
given $\underline{m}\in (0,\hat{m})$, there exist a positive
constant $\delta=\delta(\underline{m})>0$, such that if
$m\in[\underline{m},\hat{m})$, then
\begin{equation*}
q(m)=\inf_{\Omega}|\frac{\nabla\psi}{r}|\geq \delta.
\end{equation*}
\end{enumerate}
\end{lem}
\begin{pf}
\begin{enumerate}
\item [(a)] It is essentially proved in \cite{G}, however, regularity of
solutions in our case near solid boundary is only $C^{1,\alpha}$, we
can not use Hopf Lemma for derivative of solutions.

Set $\psi=\psi_2-\psi_1$, then $\psi$ satisfies
\begin{equation*}
A_{ij}(x,r)\partial_{ij}\psi+B_i(x,r)\partial_i\psi=0,
\end{equation*}
where $A_{ij}$, $B_i\in L_{loc}^{\infty}(\Omega)$. Therefore $\psi$
achieves its maximum either on the solid boundary or at far fields.
Note that $\psi$ tends to uniform flow at far fields, thus $\psi$
achieves its maximum on the nozzle wall. Therefore,
\begin{equation*}
\frac{\partial \psi}{\partial n}\geq 0,
\end{equation*}
where $n$ is the unit outer normal of the nozzle wall. On the other
hand, $\psi_i$ achieve their maximum on the whole nozzle wall where
they are constants, so
\begin{equation*}
\frac{\partial\psi_i}{\partial n}>0, \,\,\text{and}\,\,
\frac{\partial\psi_i}{\partial \tau}=0,\,\,\text{on}\,\, T_2,
\,\,\text{for}\,\, i=1,2,
\end{equation*}
where $\tau$ is the tangential direction of $T_2$. Thus
\begin{equation*}
|\frac{\nabla\psi_2}{r}|^2\geq |\frac{\nabla\psi_1}{r}|^2.
\end{equation*}
\item[(b)] It follows from the proof of Lemma
\ref{lempositiveaxialvelocity}, for any $\underline{m}\in
(0,\hat{m})$, there exists $\delta=\delta(\underline{m})$ such that
the corresponding solution $\underline{\psi}$ with $\underline{m}$
satisfies
\begin{equation*}
\inf_{\Omega}|\frac{\partial\underline{\psi}}{\partial r}|\geq
\delta.
\end{equation*}
Furthermore, $\frac{\partial \psi}{\partial r}$ achieves its infimum
either on the nozzle wall or at far fields.

Note that on the nozzle boundary
\begin{equation*}
\frac{\partial \psi}{\partial n}=\frac{\partial \psi}{\partial r}
\sqrt{1+(f'(x))^2},
\end{equation*}
therefore, for any $m\in (\underline{m},\hat{m})$,
\begin{equation*}
|\frac{\nabla\psi}{r}|\geq \frac{\partial \psi}{\partial r} \geq
\inf\frac{\partial \psi}{\partial r} \geq \inf\frac{\partial
\underline{\psi}}{\partial r}\geq \delta
\end{equation*}
\end{enumerate}
\end{pf}

Collecting all these lemmas together, we finish the proof of
Theorem \ref{Thcritical}

\section{Subsonic-Sonic Flows}

To take limit for $m_0\rightarrow \hat{m}$, let us recall the
compensated compactness framework in \cite{XX}.
\begin{thm}\label{Compensatedcompactness}
Let $w^{\varepsilon}(x,r)=(q^{\varepsilon},
\omega^{\varepsilon})(x,r)$ be a sequence of functions satisfying
the Conditions (C):

(C.1) $0<\delta\leq q^{\varepsilon}(x,r)\leq 1$ a.e. in $\Omega$
for some positive constant $\delta$.

(C.2) $|\omega^{\varepsilon}(x,r)|\leq \hat{\omega}<
\frac{\pi}{2}$, for some constant $\hat{\omega}$ independent of
$\varepsilon$.

(C.3) $\partial_{x}\eta_{\pm}(w^{\varepsilon})
+\partial_{r}\Lambda_{\pm}(w^{\varepsilon})$ are confined in a
compact set in $H_{loc}^{-1}(\Omega)$ for the momentum
entropy-entropy flux pair
\begin{equation*}
(\eta_{+},\Lambda_{+})=(\rho q^2 \cos^2\omega+p, \rho
q^2\sin\omega\cos\omega ),\quad (\eta_{-},\Lambda_{-})=(\rho
q^2\sin\omega\cos\omega, \rho q^2\sin^2\omega+p),
\end{equation*}
where $p=p(\rho)$, and $\rho=g(q^2)$ is determined by
(\ref{densityspeed}) through Bernoulli's law.

Then there exists a subsequence $\{w^{\varepsilon_k}\}$ of
$\{w^{\varepsilon}\}$ and $w(x,r)=(q,\omega)(x,r)$ such that
\begin{eqnarray}
&(q^{\varepsilon_k},\omega^{\varepsilon_k})\rightarrow (q,\omega),\label{compvconverge}\\
&q^{\varepsilon_k}\cos\omega^{\varepsilon_k}\rightarrow
q\cos\omega,\quad
q^{\varepsilon_k}\sin\omega^{\varepsilon_k}\rightarrow
q\sin\omega,\label{compirroconverge}\\
&g((q^{\varepsilon_k})^2)q^{\varepsilon_k}\cos\omega^{\varepsilon_k}\rightarrow
g(q^2)q\cos\omega,\quad
g((q^{\varepsilon_k})^2)q^{\varepsilon_k}\sin\omega^{\varepsilon_k}\rightarrow
g(q^2)q\sin\omega,\label{compmomconverge}
\end{eqnarray}
where all the convergence in
(\ref{compvconverge})-(\ref{compmomconverge}) are almost
everywhere convergence, and $w=(q,\omega)$ satisfies
\begin{eqnarray*}
0<\delta \leq q(x,r)&\leq& 1,\\
|\omega(x,r)|&\leq &\hat{\omega}.
\end{eqnarray*}
\end{thm}

\begin{rmk}
{\rm The strong convergence of velocity field $(U,V)=(q\cos
\omega,q\sin \omega)$ instead of $(q,\omega)$ was first proved in
\cite{CDSW}. Since we have good control on flow speed, we can also
get strong convergence on flow angles. The difference between
assumptions on Theorem \ref{Compensatedcompactness} and Theorem 1
in \cite{CDSW} is that they use one more entropy-entropy flux pair
instead of the condition on the lower bound on flow speed.}
\end{rmk}

Let $\psi$ satisfies (\ref{problem}). Set, as before,
\begin{equation*}
U=\frac{\psi_r}{r},\,\,
V=-\frac{\psi_x}{r},\,\,q^2=U^2+V^2,\,\,\rho=g(q^2).
\end{equation*}
By direct calculations, it is easy to see that $(\rho, U, V)$
satisfies
\begin{equation*}
\left\{
\begin{array}{ll}
(\rho U^2+p(\rho))_x+(\rho U V)_r =-\frac{\rho U V}{r},\\
(\rho U V)_x+(\rho V^2+p(\rho))_r=-\frac{\rho V^2}{r}.
\end{array}
\right.
\end{equation*}

Let $m_0^{\varepsilon}\rightarrow \hat{m}$, and
$(q^{\varepsilon},\omega^{\varepsilon})$ be the solutions to
(\ref{problem}) corresponding to $m_0^{\varepsilon}$, then away
from the symmetry axis, on any compact subset in $\Omega$,
$(-\frac{g((q^{\varepsilon})^2) U^{\varepsilon}
V^{\varepsilon}}{r}, -\frac{g((q^{\varepsilon})^2)
(V^{\varepsilon})^2}{r})$ are uniformly bounded, therefore,
precompact in $H^{-1}_{loc}$. Furthermore, due to Lemma
\ref{lemangleestimate}, the associated flow angles satisfy the
estimates
\begin{equation*}
-\frac{\pi}{2}<\min(\inf_x \arctan f'(x),0)\leq
\omega^{\varepsilon}\leq \max(\sup_x \arctan
f'(x),0)<\frac{\pi}{2}.
\end{equation*}
Therefore, conditions (C.1), (C.2) and (C.3) in Theorem
\ref{Compensatedcompactness} are all satisfied. Hence it follows
from Theorem \ref{Compensatedcompactness}, that
$(U^{\varepsilon},V^{\varepsilon})$ has weak-$\ast$ limit $(U,V)$
such that
\begin{equation}\label{limitaxialsyflow}
\left\{
\begin{array}{ll}
(r g(q^2)U)_x+(r g(q^2)V)_r=0,\\
U_r-V_x=0
\end{array}
\right.
\end{equation}
holds in $\Omega$ in the sense of distribution, where
$q^2=U^2+V^2$, and $g$ is defined in (\ref{densityspeed}).

We now verify that $(\rho, U, V)$ gives rise to a global
subsonic-sonic weak solution to (\ref{3Dproblem}) on $D$. Thus,
define
\begin{equation*}
\rho(x,y,z)=g(q^2)(x,r),\,\,u(x,y,z)=U(x,r), \,\,
v(x,y,z)=V(x,r)\frac{y}{r},\,\, w(x,y,z)=V(x,r)\frac{z}{r},
\end{equation*}
where $r=\sqrt{y^2+z^2}$. First, note that for $\eta\in
C_0^{\infty}(D_0)$,
\begin{eqnarray*}
&&\iiint_{D}u\eta_y-v\eta_x dxdydz\\
&=&\int_0^{2\pi}\int_{-\infty}^{\infty}
\int_0^{f(x)}\left(U(x,r)(\eta_r(x,r,\theta)
\cos\theta-\eta_{\theta}(x,r,\theta)\frac{\sin\theta}{r})\right)rdrdxd\theta\\
&&-\int_0^{2\pi}\int_{-\infty}^{\infty} \int_0^{f(x)}
V(x,r)\cos\theta\eta_x(x,r,\theta)rdrdxd\theta\\
&=&\int_{-\infty}^{\infty}\int_0^{f(x)}\int_0^{2\pi}
r(U(x,r)\eta_r\cos\theta-V\eta_x\cos\theta)-U\eta_{\theta}\sin\theta
d\theta drdx\\
&=&\int_{-\infty}^{\infty}\int_0^{f(x)}\left(U(x,r)\frac{\partial}{\partial
r}\int_0^{2\pi}r\eta(x,r,\theta)\cos\theta d\theta-
V(x,r)\frac{\partial}{\partial
x}\int_0^{2\pi}r\eta(x,r,\theta)\cos\theta d\theta\right)dxdr\\
&=&0,
\end{eqnarray*}
where we have used the facts that
\begin{equation*}
\int_0^{2\pi}r\eta(x,r,\theta)\cos\theta d\theta \in
C_{0}^{\infty}(\Omega)
\end{equation*}
and (\ref{limitaxialsyflow}) holds in the sense of distribution in
$\Omega$.

Define $\zeta=\zeta(s)$ such that
\begin{equation*}
\zeta(s)=1\,\,\text{if}\,\,s>1, \,\,\zeta(s)=0,\,\,
\text{if}\,\,s<\frac{1}{4}, \,\,|\zeta'|<2,\,\,\text{and}\,\,
\zeta\in C^{\infty}(\mathbb{R}),
\end{equation*}
and set $\zeta_{\delta}(s)=\zeta(\frac{s}{\delta})$,  where
$0<\delta<b$.

Then for any $\eta\in C_{0}^{\infty}(D)$,
$\zeta_{\delta}(r)\eta\in C_0^{\infty}(D_0)$, thus
\begin{eqnarray*}
&&\iiint_{D}u\eta_y-v\eta_x dxdydz\\
&=&\iint_{D}u(\zeta_{\delta}(r)\eta)_y-v(\zeta_{\delta}(r)\eta)_xdxdydz \\
&&+\iint_{D}u((1-\zeta_{\delta}(r))\eta)_y-v((1-\zeta_{\delta}(r))\eta)_xdxdydz \\
&=&\iint_{D}u((1-\zeta_{\delta}(r))\eta)_y-v((1-\zeta_{\delta}(r))\eta)_xdxdydz \\
&=&\int_0^{2\pi}d\theta\int_{-\infty}^{\infty}\int_{0}^{\delta}
\left(U(1-\zeta_{\delta}(r))_y\eta+U(1-\zeta_{\delta}(r))\eta_y-V\cos\theta
(1-\zeta_{\delta}(r))\eta_x\right)rdxdr\\
&\rightarrow& 0
\end{eqnarray*}
as $\delta\rightarrow 0$. Thus for any $\eta\in C_{0}^{\infty}(D)$
\begin{equation*}
\iiint_{D}u\eta_y-v\eta_x dxdydz=0.
\end{equation*}
Similarly, one can show that
\begin{equation*}
\iiint_{D}u\eta_z-w\eta_x dxdydz=0,
\end{equation*}
and
\begin{equation*}
\iiint_{D}v\eta_z-w\eta_y dxdydz=0,
\end{equation*}
for any $\eta\in C_{0}^{\infty}(D)$.

Finally, we check the continuity equation. For any $\eta\in
C_c^{\infty}(\bar{D})$,
\begin{eqnarray*}
&&\iiint_{D} \rho u\eta_x+\rho v\eta_y +\rho w\eta_z dxdydz\\
&=& \iiint_{D} \rho u(\zeta_{\delta}(r)\eta)_x+\rho
v(\zeta_{\delta}(r)\eta)_y
+\rho w(\zeta_{\delta}(r)\eta)_z dxdydz\\
&&+\iiint_{D} \rho u((1-\zeta_{\delta}(r))\eta)_x+\rho
v((1-\zeta_{\delta}(r))\eta)_y +\rho
w((1-\zeta_{\delta}(r))\eta)_z dxdydz
\end{eqnarray*}
Note that (\ref{problem}) shows that
\begin{equation*}
\iiint_{D}\rho^{\varepsilon}u^{\varepsilon}(\zeta_{\delta}\eta)_x
+\rho^{\varepsilon}v^{\varepsilon}(\zeta_{\delta}\eta)_y
+\rho^{\varepsilon}w^{\varepsilon}(\zeta_{\delta}\eta)_z dxdydz=0,
\end{equation*}
while (\ref{compmomconverge}) yields
\begin{eqnarray*}
&&\lim_{\varepsilon\rightarrow 0} \iiint_{D}\zeta_{\delta}(f(x)-r)
\left(\rho^{\varepsilon}u^{\varepsilon}(\zeta_{\delta}\eta)_x
+\rho^{\varepsilon}v^{\varepsilon}(\zeta_{\delta}\eta)_y
+\rho^{\varepsilon}w^{\varepsilon}(\zeta_{\delta}\eta)_z \right)dxdydz\\
&=&\iiint_{D}\zeta_{\delta}(f(x)-r)\left( \rho
u(\zeta_{\delta}(r)\eta)_x+\rho v(\zeta_{\delta}(r)\eta)_y +\rho
w(\zeta_{\delta}(r)\eta)_z\right) dxdydz,
\end{eqnarray*}
where
$(\rho^{\varepsilon},u^{\varepsilon},v^{\varepsilon},w^{\varepsilon})$
denotes the three dimensional flow associated with the incoming
mass flux $m^{\varepsilon}$. Therefore,
\begin{eqnarray*}
&& \iiint_{D} \rho u(\zeta_{\delta}(r)\eta)_x+\rho
v(\zeta_{\delta}(r)\eta)_y +\rho w(\zeta_{\delta}(r)\eta)_z dxdydz\\
&=&\iiint_{D}(1-\zeta_{\delta}(f(x)-r))( \rho
u(\zeta_{\delta}(r)\eta)_x+\rho
v(\zeta_{\delta}(r)\eta)_y +\rho w(\zeta_{\delta}(r))\eta)_z dxdydz\\
&&+\iiint_{D}\zeta_{\delta}(f(x)-r)( \rho
u(\zeta_{\delta}(r)\eta)_x+\rho
v(\zeta_{\delta}(r)\eta)_y +\rho w(\zeta_{\delta}(r))\eta)_z dxdydz\\
&=&\iiint_{D}(1-\zeta_{\delta}(f(x)-r))( \rho
u(\zeta_{\delta}(r)\eta)_x+\rho
v(\zeta_{\delta}(r)\eta)_y +\rho w(\zeta_{\delta}(r))\eta)_z dxdydz\\
&&+\lim_{\varepsilon\rightarrow 0} \iiint_{D}\zeta_{\delta}(f(x)-r)
\left(\rho^{\varepsilon}u^{\varepsilon}(\zeta_{\delta}\eta)_x
+\rho^{\varepsilon}v^{\varepsilon}(\zeta_{\delta}\eta)_y
+\rho^{\varepsilon}w^{\varepsilon}(\zeta_{\delta}\eta)_z \right)
dxdydz
\end{eqnarray*}
\begin{eqnarray*}
&=&\iiint_{D}(1-\zeta_{\delta}(f(x)-r))( \rho
u(\zeta_{\delta}(r)\eta)_x+\rho
v(\zeta_{\delta}(r)\eta)_y +\rho w(\zeta_{\delta}(r))\eta)_z dxdydz\\
&&+\lim_{\varepsilon\rightarrow 0} \iiint_{D}
\left(\rho^{\varepsilon}u^{\varepsilon}(\zeta_{\delta}\eta)_x
+\rho^{\varepsilon}v^{\varepsilon}(\zeta_{\delta}\eta)_y
+\rho^{\varepsilon}w^{\varepsilon}(\zeta_{\delta}\eta)_z \right) dxdydz\\
&&+\lim_{\varepsilon\rightarrow
0}\iiint_{D}(\zeta_{\delta}(f(x)-r)-1)
\left(\rho^{\varepsilon}u^{\varepsilon}(\zeta_{\delta}\eta)_x
+\rho^{\varepsilon}v^{\varepsilon}(\zeta_{\delta}\eta)_y
+\rho^{\varepsilon}w^{\varepsilon}(\zeta_{\delta}\eta)_z \right)
dxdydz\\
 &=&\iiint_{D}(1-\zeta_{\delta}(f(x)-r))( \rho
u(\zeta_{\delta}(r)\eta)_x+\rho
v(\zeta_{\delta}(r)\eta)_y +\rho w(\zeta_{\delta}(r))\eta)_z dxdydz\\
&&+\lim_{\varepsilon\rightarrow
0}\iiint_{D}(\zeta_{\delta}(f(x)-r)-1)
\left(\rho^{\varepsilon}u^{\varepsilon}(\zeta_{\delta}\eta)_x
+\rho^{\varepsilon}v^{\varepsilon}(\zeta_{\delta}\eta)_y
+\rho^{\varepsilon}w^{\varepsilon}(\zeta_{\delta}\eta)_z \right)
dxdydz.
\end{eqnarray*}
Thus
\begin{eqnarray*}
&&\iiint_{D} \rho u\eta_x+\rho v\eta_y +\rho w\eta_z dxdydz\\
&=&\iiint_{D} \rho u((1-\zeta_{\delta}(r))\eta)_x+\rho
v((1-\zeta_{\delta}(r))\eta)_y +\rho
w((1-\zeta_{\delta}(r))\eta)_z dxdydz\\
&&+\iiint_{D}(1-\zeta_{\delta}(f(x)-r))\left( \rho
u(\zeta_{\delta}(r)\eta)_x+\rho
v(\zeta_{\delta}(r)\eta)_y +\rho w(\zeta_{\delta}(r)\eta)_z\right) dxdydz\\
&&+\lim_{\varepsilon\rightarrow
0}\iiint_{D}(\zeta_{\delta}(f(x)-r)-1)
\left(\rho^{\varepsilon}u^{\varepsilon}(\zeta_{\delta}\eta)_x
+\rho^{\varepsilon}v^{\varepsilon}(\zeta_{\delta}\eta)_y
+\rho^{\varepsilon}w^{\varepsilon}(\zeta_{\delta}\eta)_z \right)\\
&\rightarrow& 0
\end{eqnarray*}
as $\delta\rightarrow 0$. Therefore, for any test function
$\eta\in C_c^{\infty}(\bar{D})$,
\begin{equation}\label{limit3Dpotentialeq}
\iiint_{D} \rho u\eta_x+\rho v\eta_y +\rho w\eta_z dxdydz=0.
\end{equation}
Moreover, the equation (\ref{limit3Dpotentialeq}) also implies
that $(u, v, w)$ satisfies the boundary condition (\ref{3Dbc})
actually as the normal trace of divergence measure field $(\rho u,
\rho v, \rho w)$ on the boundary in the sense of
Anzellotti\cite{Anzellotti}. So we finish the proof of Theorem 3.

Further characterizations of the subsonic-sonic flow we obtained
are left for future.

{\bf Acknowledgment.} The authors would like to thank Mr. Wei Yan
and Professor Hongjun Yuan for many stimulating discussions. This
research is supported in part by Hong Kong RGC Earmarked Research
Grants CUHK04-28/04, CUHK-4040/06, and RGC Central Allocation
Grant CA05/06.SC01.

\end{document}